\def \NN {\mathbb N}
\def \CC {\mathbb C}

\def \RR {\mathbb R}
\def \ZZ {\mathbb Z}

\def \epsilon{\varepsilon}

\def \A  {{\mathcal A}}

\def \C  {{\mathcal C}}
\def \D  {{\mathcal D}}

\def \LL {{\mathcal L}}
\def \K {{\mathcal K}}

\def \d {\text{d}}

\def \bfz {{\boldsymbol{z}}}
\def \bfx {{\boldsymbol{x}}}
\def \bfy {{\boldsymbol{y}}}
\def \bfw {{\boldsymbol{w}}}
\def \bfT {{\boldsymbol{T}}}
\def \bfS {{\boldsymbol{S}}}
\def \bfG {{\boldsymbol{G}}}
\def \bfalpha {{\boldsymbol{\alpha}}}

\def \fine {{\hfill \qedsymbol}}

\def \si {\sigma}

\newcommand{\res}{\text{res}}
\renewcommand{\S}{{\mathcal S}}

\documentclass[12pt,reqno]{amsart}

\usepackage{amsmath,amssymb}
\numberwithin{equation}{section}

\begin{document}

\title[]{Twists and resonance of $L$-functions, II}

\author[]{J.Kaczorowski  \lowercase{and} A.Perelli}
\date{}
\maketitle

\bigskip
{\bf Abstract.} We continue our investigations of the analytic properties of nonlinear twists of $L$-functions developed in \cite{Ka-Pe/2005},\cite{Ka-Pe/2011a} and \cite{Ka-Pe/resoI}. Let $F(s)$ be an $L$-function of degree $d$. First we extend the transformation formula in \cite{Ka-Pe/2011a}, relating a twist $F(s;f)$ with leading exponent $\kappa_0>1/d$ to its dual twist $\overline{F}(s;f^*)$. Then we combine the results in  \cite{Ka-Pe/resoI} with such a transformation formula to obtain the analytic properties of new classes of nonlinear twists. This allows to detect several new cases of resonance of the classical $L$-functions.

\medskip
{\bf Mathematics Subject Classification (2000):} 11 M 41

\medskip
{\bf Keywords:} $L$-functions, Selberg class, twists.

\vskip1cm
\section{Introduction}

\smallskip
1.1.~{\bf Outline.} This paper is a continuation of \cite{Ka-Pe/resoI}, to which we refer for definitions, notation and a general discussion of twists of $L$-functions. In \cite{Ka-Pe/resoI} we gave a rather complete description of the analytic properties of the nonlinear twists (as usual $e(x)=e^{2\pi ix}$)
\[
F(s;f) = \sum_{n=1}^\infty \frac{a(n)}{n^s} e(-f(n,\bfalpha))
\]
with leading exponent $\kappa_0\leq 1/d$, as well as some applications to the resonance problem for such twists. Here $F(s)$ is a function of degree $d$ in the extended Selberg class (i.e. $F\in\S_d^\sharp$) and 
\begin{equation}
\label{1-1}
f(\xi,\bfalpha) = \sum_{j=0}^N \alpha_j \xi^{\kappa_j}, \hskip1.5cm 0\leq \kappa_N<\dots<\kappa_0, \  \alpha_j\in\RR \ \text{and} \ \prod_{j=0}^N\alpha_j\neq0. 
\end{equation}
In this paper we deal with the case $\kappa_0>1/d$. We start with an extension of the transformation formula, given in Theorem 1.1 of \cite{Ka-Pe/2011a}, relating a nonlinear twist $F(s;f)$ with $\kappa_0>1/d$ to its dual twist $\overline{F}(s;f^*)$; see Theorem 1 below. Then we combine the results in \cite{Ka-Pe/resoI} and Theorem 1 to obtain the analytic properties of new classes of nonlinear twists. Finally, we consider in greater detail some classes of twists of degree 2 $L$-functions. The applications to the resonance problem are similar to those in \cite{Ka-Pe/resoI}, hence we shall very briefly outline them.

\medskip
In order to give the flavor of the results in the paper, we first present two very special cases of our theorems. The elliptic curve $E$ of equation $y^2-y = x^3-x$ has conductor 37 and corresponds to a newform of weight 2 and level 37; see Zagier \cite{Zag/1985} for a pre-Wiles justification of this. Denoting by $L_E(s)$ the associated $L$-function and by $a_E(n)$ its coefficients, we have that
\[
F(s) = L_E(s+\frac12)
\]
is a degree 2 element of the Selberg class, hence we may apply Theorem 6 below. Choosing $k=1$, $\ell=0$ and $\alpha=(37)^{-2/3}$, and recalling that the conductor of $F(s)$ equals 37, from the definition of spectrum of $F(s)$ in \cite{Ka-Pe/resoI} we see that $2(37)^{-1/2} \in \ \text{Spec}(F)$. Hence, using the notation introduced below, from Theorem 6 we deduce that the corresponding function 
\[
f(\xi)= \frac{3}{2(2738)^{1/3}} \xi^{2/3} +  \frac{1}{(1369)^{1/3}} \xi^{1/6}
\]
belongs to $\A_0(F)\setminus \A_{00}(F)$, since $F(s)$ is entire. Thus by Theorem 4 below (after a change of variable) we see that the twisted $L$-function
\[
L_E(s;f)=\sum_{n=1}^\infty \frac{a_E(n)}{n^s} e(-\frac{3}{2(2738)^{1/3}} n^{2/3} - \frac{1}{(1369)^{1/3} }n^{1/6})
\]
is meromorphic on $\CC$, and has a simple pole at $s_0 = 1 + \frac{1}{4D(f)}$ since $\theta_F=0$. Moreover, we have that $f=TS(f_0)$, where $f_0$ is the standard twist and $S$ is a shift of degree 2; hence $D(f) = 2\cdot2-1=3$, thus $s_0=13/12$. As a consequence, choosing $r=1$ in \eqref{1-16} below, we deduce the following asymptotic formula for the associated nonlinear exponential sum:
\[
\sum_{n=1}^\infty a_E(n) e(-\frac{3}{2(2738)^{1/3}} n^{2/3} - \frac{1}{(1369)^{1/3} }n^{1/6}) e^{-n/x} \sim c_0 x^{13/12}
\]
as $x\to\infty$, with a certain $c_0\neq0$. If we change the value of $\alpha$, e.g. we choose $\alpha=(37)^{-20/31}$, then by a similar argument we see that the resulting function $f(\xi)$ belongs to $\A(F)\setminus\A_0(F)$. Hence by Theorems 3 and 6 below (again after a change of variable) we obtain that
\[
L_E(s;f)=\sum_{n=1}^\infty \frac{a_E(n)}{n^s} e(-\frac{3}{2(2738)^{1/3}} n^{2/3} - \frac{1}{(37)^{20/31} }n^{1/6})
\]
is entire, and the corresponding nonlinear exponential sum is $O(1)$. Similar results can be deduced, from the above quoted theorems, for the coefficients any degree 2 $L$-function, e.g. the divisor function $d(n)$ or the Ramanujan function $\tau(n)$.

\medskip
1.2.~{\bf Transformation formula.} In order to state Theorem 1 we need to introduce further notation, not present in \cite{Ka-Pe/resoI}, concerning nonlinear twists of a given $F\in\S_d^\sharp$, as well as some results proved in \cite{Ka-Pe/2011a}. First we rewrite the twist function $f(\xi,\bfalpha)$ (sometimes simply $f(\xi)$ or $f$) as
\[
f(\xi,\bfalpha) = \xi^{\kappa_0}\sum_{j=0}^N \alpha_j \xi^{-\omega_j} \hskip1.5cm 0=\omega_0<\dots<\omega_N\leq \kappa_0,
\]
we recall that $\kappa_0$ is the {\it leading exponent} of $f$ and write $\kappa_0=\text{lexp}(f)$, and consider the semigroup
\[
\D_f = \{\omega=\sum_{j=1}^N m_j\omega_j: m_j\in\ZZ, \ m_j\geq 0\}.
\]
If $\alpha_0>0$ and $ \kappa_0>1/d$ we put ($q_F$ is the conductor of $F(s)$, see e.g. p.1398 of \cite{Ka-Pe/2011a})
\begin{equation}
\label{1-2}
\Phi(z,\xi,\bfalpha) = z^{1/d} - 2\pi f(\frac{q z}{\xi},\bfalpha) \hskip1.5cm q=q_F(2\pi d)^{-d}.
\end{equation}
Thanks to (i) of Theorem 1.2 of \cite{Ka-Pe/2011a} we have the operator
\begin{equation}
\label{1-3}
T(f)(\xi,\bfalpha) := \frac{1}{4\pi^2 i} \int_{\C} \frac{\Phi(z,\xi,\bfalpha) \frac{\partial^2}{\partial z^2} \Phi(z,\xi,\bfalpha)}{\frac{\partial}{\partial z} \Phi(z,\xi,\bfalpha)} \d z = \xi^{\kappa_0^*} \sum_{\omega\in\D_f} A_\omega(\bfalpha) \xi^{-\omega^*},
\end{equation}
where
\[
\kappa_0^* = \frac{\kappa_0}{d\kappa_0-1}, \hskip1.5cm \omega^* = \frac{\omega}{d\kappa_0-1}
\]
and $\C$ is, roughly, a circle around the $z$-critical point of $\Phi(z,\xi,\bfalpha)$; see p.1401 of \cite{Ka-Pe/2011a} for the exact definition of $\C$. Moreover, from Theorem 1.2 of \cite{Ka-Pe/2011a} we obtain some properties of the coefficients $A_\omega(\bfalpha)$, and also that $T$ is self-reciprocal (i.e. $T^2=$ identity).

\medskip
{\bf Remark 1.} Note that the operator $T$ depends on the degree $d$ of $F(s)$ but also on its conductor $q_F$. Therefore we fix the function $F\in\S^\sharp_d$ when dealing with the operator $T$ or, at least, we fix degree and conductor. \fine

\medskip
\noindent
We also write
\begin{equation}
\label{1-4}
s^* = \frac{s+\frac{d\kappa_0}{2}-1+id\theta_F\kappa_0}{d\kappa_0-1},
\end{equation}
where $\theta_F$ is the internal shift of $F(s)$ (see \cite{Ka-Pe/resoI}; usually $\theta_F=0$ for the classical $L$-functions). We denote by $T^\flat(f)(\xi,\bfalpha)$ the finite sum obtained from the right hand side of \eqref{1-3} dropping the terms with $\omega>\kappa_0$, and write
\begin{equation}
\label{1-5}
f^*(\xi,\bfalpha) = T^\flat(f)(\xi,\bfalpha).
\end{equation}
Since $T$ is self-reciprocal we have
\[
(f^*)^*(\xi,\bfalpha) = f(\xi,\bfalpha).
\]
Finally, for $\alpha_0<0$ we put 
\[
T(f)(\xi,\bfalpha) = -T(-f)(\xi,\bfalpha),
\]
and hence for every $f(\xi,\bfalpha)$ we have
\begin{equation}
\label{1-6}
f^*(\xi,\bfalpha) = - (-f)^*(\xi,\bfalpha) \hskip1.5cm (f^*)^*(\xi,\bfalpha) = f(\xi,\bfalpha).
\end{equation}

\medskip
{\bf Theorem 1.} {\sl Let $F\in\S_d^\sharp$ with $d\geq1$, $f(\xi,\bfalpha)$ be as in \eqref{1-1} with $\kappa_0>1/d$, and let $K\geq0$. There exist an integer $J\geq0$, constants $0=\eta_0<\eta_1<\dots<\eta_J$ and functions $W_0(s),\dots,W_J(s)$, $G(s;f)$ holomorphic for $\sigma>-K$, with $W_0(s)$ nonvanishing there, such that for $\sigma>-K$
\begin{equation}
\label{1-7}
F(s;f) = \sum_{j=0}^J W_j(s) \overline{F}(s^*+\eta_j;f^*) + G(s;f).
\end{equation}
Moreover, uniformly for $s$ in any given vertical strip inside $\sigma>-K$ and any $\epsilon>0$ we have}
\[
G(s;f) \ll e^{\epsilon|t|} \hskip1.5cm |t|\to\infty.
\]

\medskip
The meaning of \eqref{1-7} is that $F(s;f) - \sum_{j=0}^J W_j(s) \overline{F}(s^*+\eta_j;f^*)$ has holomorphic continuation to $\si>-K$. If $1<d\kappa_0<2$ and $\si>1/2$, then $\Re s^* > \si$ and hence \eqref{1-7} gives some analytic continuation of $F(s;f)$ to the left of $\si=1$. However, below we show that Theorem 1 can be coupled with other ideas to give meromorphic continuation for a new class of nonlinear twists.

\medskip
{\bf Remark 2.}  Comparing with Theorem 1.1 of \cite{Ka-Pe/2011a}, we see that in Theorem 1 we dropped the three restrictions $F(s)$ entire, $\theta_F=0$ and $\sigma>0$. Hence the transformation formula in Theorem 1.1 of \cite{Ka-Pe/2011a} holds now for every $F\in\S_d^\sharp$ and for $s$ in any given right half-plane. Moreover, Theorem 1 contains in addition a bound on the order of growth of the function $G(s;f)$. Note that condition $d\geq 1$ is equivalent to $d>0$, since the extended Selberg class is empty for degrees between 0 and 1; see Conrey-Ghosh \cite{Co-Gh/1993} and our paper \cite{Ka-Pe/1999a}. Note also that the integer $J$, the constants $\eta_j$ and the functions $W_j(s)$ and $G(s;f)$ may depend on $F(s)$, $f(\xi,\bfalpha)$ and $K$.  \fine

\medskip
{\bf Remark 3.} An inspection of the proof of Theorem 1, see Section 2.5, allows to say something more on the shifts $\eta_j$ and on the functions $W_j(s)$. In particular, $\eta_J>(K+\frac{d\kappa_0}{2})/(d\kappa_0-1)$ and the  $W_j(s)$ are of the form $W_j(s) = e^{a_js}P_j(s)$ with $a_j\in\RR$, $P_j\in\CC[s]$ and deg$P_j \ll j$. \fine

\medskip
{\bf Remark 4.} In Theorem 4 of \cite{Ka-Pe/resoI} we showed that the presence of negative exponents in the function $f(\xi,\bfalpha)$ produces a kind of stratification of the twist $F(s;f)$ in terms of shifts of the twists $F(; f^\flat)$ ($f^\flat$ denotes, according to the above notation, the part of $f$ with exponents $\geq0$); see formula (4.1) there. In the proof of Theorem 1 we actually cut off the negative exponents coming out naturally in $T(f)(\xi,\bfalpha)$, but their effect is present in the stratification appearing in formula \eqref{1-7} above. \fine

\medskip
The strategy of the proof of Theorem 1 is sketched at the beginning of Section 2.

\medskip
1.3.~{\bf The group $\frak{G}$.} Now we turn to the applications of Theorem 1 and of the results in \cite{Ka-Pe/resoI}, but we first have to set up further notation. Along with the {\it duality} operator $T$ in \eqref{1-3}, whose applicability is limited by the fact that it is self-reciprocal, we consider the {\it shift} operator $S$, defined by
\begin{equation}
\label{1-8}
S(f)(\xi,\bfalpha) = f(\xi,\bfalpha) + P(\xi),
\end{equation}
where $P\in\ZZ[\xi]$ has $\deg P\geq 1$. We define the degree of $S$ as $\deg S= \deg P$, and note that every $S$ can be obtained by repeated applications of the shifts $S^{(k)}$, associated with $P(\xi)=\xi^k$ for $k=1,2,\dots$, and of their inverses (associated with $P(\xi)=-\xi^k$). Clearly, $S$ acts trivially on twists, i.e. $F(s;S(f)) = F(s;f)$, but the action of $T$ on $S(f)(\xi,\bfalpha)$ differs from the action of $T$ on $f(\xi,\bfalpha)$. This enables to build nontrivial chains of applications of the operators $T$ and $S$.

\medskip
{\bf Remark 5.} This phenomenon, along with the transformation formula (see Theorem 1) and the properties of the standard twist (see \cite{Ka-Pe/2005} and \cite{Ka-Pe/resoI}),  is the basis for the most interesting aspects of our twist theory. Indeed, in \cite{Ka-Pe/2002} and \cite{Ka-Pe/2011a} we exploited it to prove the degree conjecture for $\S^\sharp$ in the range $1<d<2$, while in this paper it is used to obtain the analytic properties of a new class of nonlinear twists. \fine

\medskip
We first consider the formal group $\frak{G}$ generated by the symbols $\bfT, \bfS^{(1)}, \bfS^{(2)}, \dots$, satisfying the relations $\bfT^2=$ identity and $\bfS^{(h)}\bfS^{(k)} = \bfS^{(k)}\bfS^{(h)}$. Thus every element of $\frak{G}$ has the form
\[
\bfG = \bfT^a \bfS_k \bfT \bfS_{k-1} \cdots \bfT \bfS_1 \bfT^b \hskip1.5cm a,b\in\{0,1\},
\]
where each $\bfS_j$ is a product of $\bfS^{(k)}$'s and their inverses. Then, for given $F\in\S_d^\sharp$ with $d\geq1$ and $f(\xi,\bfalpha)$ as in \eqref{1-1}, to an element $\bfG$ as above we associate, with obvious notation, the operator $G$ defined by
\begin{equation}
\label{1-9}
G(f) = (T^\flat)^a S_k T^\flat S_{k-1} \cdots T^\flat S_1 (T^\flat)^b (f),
\end{equation}
provided it is well defined. Indeed, in view of Remark 1 the action of $T^\flat$ depends on $F(s)$ and, by Theorem 1, $T^\flat(f)$ makes sense under the compatibility condition that lexp$(f)>1/d$. Since in the applications of Theorem 1 in this paper we usually start with functions $f_0(\xi,\bfalpha)$ satisfying
\begin{equation}
\label{1-10}
0\leq \text{lexp}(f_0) \leq 1/d,
\end{equation}
we choose $b=0$ in \eqref{1-9}. Moreover, we may always assume that $a=1$ since, as we said, the action of $S$ on twists is trivial. Therefore we may assume that $G$ has the form
\begin{equation}
\label{1-11}
G = T^\flat S_k T^\flat S_{k-1} \cdots T^\flat S_1,
\end{equation}
in which case $G(f_0)$ is well defined if and only if for $j=1,\dots,k$ we have
\begin{equation}
\label{1-12}
\ell_j := \text{lexp}(S_jT^\flat \cdots T^\flat S_1(f_0))> 1/d.
\end{equation}
We remark that condition \eqref{1-12} is not empty only when $d\leq 2$, i.e. for $d=1$ and $d=2$ in view of the main result in \cite{Ka-Pe/2011a}, thanks to the following   

\smallskip
{\bf Fact 1.} {\sl If $d>2$ then $\ell_j=\deg S_j$ for $j=1,\dots,k$. In particular, \eqref{1-12} is satisfied.} \fine

\smallskip
\noindent
The proof of this fact will be given at the beginning of Section 3. 

\medskip
Given $F\in\S^\sharp$, if $G$ as in \eqref{1-11} and $f_0(\xi,\bfalpha)$ as in \eqref{1-1} satisfy \eqref{1-10} and \eqref{1-12}, thanks to Theorem 1 and the results in \cite{Ka-Pe/resoI} we may hope to get nontrivial information of the analytic properties of the twist $F(s;G(f_0))$. We therefore define
\begin{equation}
\label{1-13}
\A(F) = \{f = G(f_0): \text{$G$ as in \eqref{1-11} and $f_0$ as in \eqref{1-10} satisfy \eqref{1-12}}\}.
\end{equation}
Given $G(f_0), H(f_1)\in\A(F)$, we first observe the following

\smallskip
{\bf Fact 2.} {\sl If $d>2$ we have $G(f_0)=H(f_1) \Rightarrow G=H$.} \fine

\smallskip
\noindent
Again, for the proof see the beginning of Section 3. For $d\leq 2$ this does not hold, as the following example shows. Suppose that $F(s)$ has degree 2 and conductor 1, and let $f_0(\xi)=\alpha\sqrt{\xi}$. Then a computation based on Theorem 1.2 of \cite{Ka-Pe/2011a} shows that 
\[
T^\flat S^{(1)}(f_0)(\xi) = -\xi +\alpha\sqrt{\xi},
\]
hence, for example, $T^\flat S^{(1)} S^{(1)}T^\flat S^{(1)}(f_0) = T^\flat S^{(1)}(f_0)$ but $T^\flat S^{(1)} S^{(1)}T^\flat S^{(1)} \neq T^\flat S^{(1)}$.

\medskip
The number of $T^\flat$ in $G$ is called the {\it weight} of $G$ and is denoted by $\omega(G)$; moreover, given $f\in\A(F)$ we denote by
\[
\omega(f) = \min\{\omega(G):  \text{$G$ as in \eqref{1-11} and $G(f_0)=f$ for some $f_0$ satisfying \eqref{1-10}}\}
\]
the {\it weight} of $f$. In view of Fact 2, for $d>2$ the pair $(G,f_0)$ is uniquely determined by $f=G(f_0)$, and in this case we have $\omega(f)=\omega(G)$. If $f\in\A(F)$ has the form $f=G(f_0)$ and $\omega(f)\geq 1$ we write
\[
D(f) = \prod_{j=1}^k(d\ell_j-1),
\]
where the $\ell_j$ are associated with $f_0$ as in \eqref{1-12}, while if $\omega(f)=0$ we simply put $D(f) = 1$. The main feature of $D(f)$ is that it is independent of the particular representation of $f$ inside $\A(F)$, as shown (with obvious notation) by the following

\smallskip
{\bf Fact 3.} {\sl If $f=G(f_0)=H(f_0')$ then $\prod_{j=1}^k (d\ell_j-1) = \prod_{j=1}^{k'} (d\ell'_j-1)$.} \fine

\smallskip
\noindent
Fact 3 follows at once from Theorem 4 below. Of course, in view of Fact 2, it is nontrivial only for $d\leq 2$. 

\medskip
Finally, we refer to \cite{Ka-Pe/resoI} for the notions and properties of standard twist $F(s,\alpha)$ and spectrum Spec$(F)$ of $F\in\S^\sharp$, and define
\[
\text{Spec$^*(F) = \{\alpha\in\RR: F(s,\alpha)$ is not entire$\}$}.
\]
Clearly
\[
\text{Spec}^*(F) = 
\begin{cases}
\text{Spec}(F) \cup (-\text{Spec}(F)) \ \ &  \text{if $F(s)$ is entire} \\
\text{Spec}(F) \cup (-\text{Spec}(F)) \cup\{0\} \ \ & \text{if $F(s)$ has a pole at $s=1$.}
\end{cases}
\]
The results in \cite{Ka-Pe/resoI} show a different behavior of $F(s,\alpha)$ with $\alpha\in$ Spec$(F)$ with respect to the other twists $F(s;f_0)$ with $0\leq \text{lexp}(f_0) \leq 1/d$; accordingly, we consider the subsets of $\A(F)$, see \eqref{1-13}, defined as
\begin{equation}
\label{1-14}
\A_0(F) = \{f\in\A(F): f=G(f_0) \ \text{with} \ f_0(\xi,\alpha) = \alpha\xi^{1/d}, \alpha\in \text{Spec}^*(F)\}
\end{equation}
\begin{equation}
\label{1-15}
\A_{00}(F) = \{f\in\A_0(F) \ \text{with $\alpha=0$ in \eqref{1-14}}\}.
\end{equation}
Of course $\A_{00}(F) \subset \A_0(F) \subset \A(F)$, and $\A_{00}(F) \neq \emptyset$ if and only if $F(s)$ is polar; recall that the polar order of $F(s)$ at $s=1$ is denoted by $m_F$. 

\medskip
1.4.~{\bf Nonlinear twists.} Now we can state several general results about the twists $F(s;f)$. Recalling the definitions in \eqref{1-13}, \eqref{1-14} and \eqref{1-15} we have

\medskip
{\bf Theorem 2.} {\sl Let $F\in\S^\sharp_d$ with $d\geq 1$. Then $F(s;f)$ is meromorphic on $\CC$ for every $f\in\A(F)$, with all poles in a horizontal strip of type $|t|\leq T_0(F)$ and of order $\leq \max(1,m_F)$. Moreover, for any $\epsilon>0$ and $A<B$, as $|t|\to\infty$ we have
\[
F(s;f) \ll e^{\epsilon|t|}
\]
uniformly for} $A\leq \sigma\leq B$.

\medskip
{\bf Theorem 3.} {\sl Let $F\in\S^\sharp_d$ with $d\geq 1$. If $f\in\A(F)\setminus \A_0(F)$ then $F(s;f)$ is entire.}

\medskip
Probably, in this case the twists $F(s;f)$ have finite order. Proving this would require additional uniformity estimates in Theorem 1, similar to those in \cite{Ka-Pe/twist}. Since this would considerably enlarge the size of the paper, we shall not enter such a problem.

\medskip
{\bf Theorem 4.} {\sl Let $F\in\S^\sharp_d$ with $d\geq 1$. If $f\in\A_0(F)\setminus \A_{00}(F)$ then $F(s;f)$ is not entire, all its poles are simple and lie on the half-line 
\[
s=\sigma -i\theta_F \quad \text{with} \quad  \sigma \leq \frac12 + \frac{1}{2dD(f)},
\]
the initial point $s_0= \frac12 + \frac{1}{2dD(f)} -i\theta_F$ being a pole.}

\medskip
{\bf Theorem 5.} {\sl Let $F\in\S^\sharp_d$ with $d\geq 1$. If $f\in\A_{00}(F)$ then $F(s;f)$ is not entire, all its poles have order $\leq m_F$ and lie on the half-line 
\[
s=\sigma + i\theta_F\big(\frac{(-1)^{\omega(f)}}{D(f)}-1\big)  \quad \text{with} \quad  \sigma \leq \frac12 + \frac{1}{2D(f)},
\]
the initial point $s_0= \frac12 + \frac{1}{2D(f)}  + i \big(\frac{(-1)^{\omega(f)}}{D(f)}-1\big) \theta_F$ being a pole of order $m_F$.}

\medskip
1.5.~{\bf Resonance.} As an application of the above theorems we consider the smoothed nonlinear exponential sums
\begin{equation}
\label{1-16}
S_F(x;f,r) = \sum_{n=1}^\infty a(n)e(-f(n,\bfalpha)) e^{-(n/x)^r}
\end{equation}
with $r>0$ arbitrary. Then, as $x\to\infty$, under the hypotheses of Theorems 3, 4 and 5 we have respectively that
\[
\begin{split}
S_F(x;f,r) &= O(1), \\
S_F(x;f,r) &\sim c_0 x^{1/2 + 1/(2dD(f) )- i\theta_F} 
\end{split}
\]
with a certain $c_0\neq 0$,  and
\[
\hskip1.2cm S_F(x;f,r) \sim  x^{1/2 + 1/(2D(f) )+ i\lambda_f\theta_F} P(\log x) 
\]
with $\lambda_f = \big(\frac{(-1)^{\omega(f)}}{D(f)}-1\big)$ and $P\in\CC[x]$, $\deg P = m_F-1$. Proofs are standard and hence omitted.

\medskip
1.6.~{\bf Degree 2 examples.} As an illustration of the above theorems we describe in detail some functions $f\in\A(F)$ in the degree 2 case. The two examples presented in Section 1.1, concerning the $L$-function associated with an elliptic curve, are special cases of the following

\medskip
{\bf Theorem 6.} {\sl Let $F\in\S^\sharp_2$. Then for every $k,\ell\in\ZZ$, $k\neq0$, and $\alpha\in\RR$ the function
\[
f(\xi) = \frac{3}{2} \frac{1}{(2kq_F^2)^{1/3}} \xi^{2/3} + \frac{3}{4} \frac{\ell}{(2k^2q_F)^{2/3}} \xi^{1/3} + \alpha \xi^{1/6}
\]
belongs to $\A(F)$. Moreover,} $f\in\A_0(F) \Leftrightarrow 2(k^2q_F)^{1/6}|\alpha| \in$ Spec$(F)$.

\medskip
{\bf Remark 6.}  In \cite{Ka-Pe/resoI} we proved that if the leading exponent $\kappa_0$ is $\leq 1/d$ and $f(\xi,\bfalpha)$ has at least one exponent $<1/d$, then $F(s;f)$ is entire. Theorem 6, and in particular the first example given in Section 1.1, shows that this is not true in the case $\kappa_0>1/d$. A heuristic explanation of these facts would be of interest. \fine

\medskip
We conclude remarking that explicit results similar to Theorem 6 can be computed for any $L$-function of degree $d\geq1$. For example, starting with $\zeta(s)$ we can get the basic analytic properties of
\[
\sum_{n=1}^\infty \frac{e(\alpha n^{3/2} +\beta n)}{n^s}
\]
with $\alpha$ in a certain infinite discrete set and any real $\beta$. As a consequence, the asymptotic behavior of the corresponding nonlinear exponential sums can be detected.

\medskip
{\bf Acknowledgements.} 
This research was partially supported by the Istituto Nazionale di Alta Matematica, by grant PRIN2010-11 {\sl Arithmetic Algebraic Geometry and Number Theory} and by grants N N201 605940 and 2013/11/B/ST1/ 02799 {\sl Analytic Methods in Arithmetic} of the National Science Centre.

\bigskip
\section{Proof of Theorem 1}

\smallskip
We follow the proof of Theorem 1.1 of \cite{Ka-Pe/2011a}, indicating only the main differences and referring to our previous papers \cite{Ka-Pe/2005}, \cite{Ka-Pe/2011a} and \cite{Ka-Pe/resoI} as much as possible. Moreover, we take this opportunity to streamline the arguments in \cite{Ka-Pe/2011a} and to correct some inaccuracies. Let $F(s)$ and $f(\xi,\bfalpha)$ be as in Theorem 1, with $\alpha_0>0$ otherwise we take conjugates, and for $X\geq 1$ let
\[
F_X(s;f) = \sum_{n=1}^\infty \frac{a(n)}{n^s} e(-f(n,\bfz)), \qquad \bfz = (z_0,\dots,z_N), \qquad z_j=\frac{1}{X}+2\pi i\alpha_j.
\]
Clearly, the series is convergent for every $s\in\CC$, and for $\si>1$
\[
\lim_{X\to\infty} F_X(s;f) = F(s;f).
\]

\medskip
Before embarking into the proof of Theorem 1, we briefly sketch our basic strategy. Note that the arguments in the proof sometimes give more than the minimal requirements needed to prove Theorem 1, which we use in the following sketch. Let $R$ be a given positive real number and $\epsilon>0$ be arbitrarily small. We first show, see Lemma 2.1 below, that for $-R<\sigma<-R+\delta$
\[
F_X(s;f) = M^{(1)}_X(s) + H^{(1)}_X(s),
\]
where $\delta>0$ is a small constant, $M^{(1)}_X(s)$ is a certain main term and the error term $H^{(1)}_X(s)$ is $\ll e^{\epsilon|t|}$ as $|t|\to\infty$, uniformly in $X$. Then we progressively refine $M^{(1)}_X(s)$, thus getting similar expressions for $F_X(s)$ with certain terms $M^{(2)}_X(s), \dots$ and $H^{(2)}_X(s), \dots$, satisfying the same upper bound, in place of $M^{(1)}_X(s)$ and $H^{(1)}_X(s)$, respectively. Eventually we obtain that
\[
F_X(s;f) - M_X(s;f) = H_X(s;f),
\]
where, thanks to its rather explicit form, $M_X(s;f)$ is holomorphic for $\si>-R$ and hence so is $H_X(s;f)$. Moreover, $H_X(s;f)\ll e^{\epsilon|t|}$ uniformly in $X$ for $-R<\sigma<-R+\delta$ and, again thanks to the form of $M_X(s;f)$, $H_X(s;f)$ is bounded uniformly in $X$ for $\si>\si_0$, with a certain $\si_0>1$. Further, for each $X\geq1$ we have that $H_X(s;f)$ is bounded for $\si>-R$. Therefore, by an application of the Phragm\'en-Lindel\"of theorem we have that $H_X(s;f)\ll e^{\epsilon|t|}$ uniformly in $X$ for $s$ in any substrip of $-R<\sigma<\si_0+1$. Since for $\si>\si_0$ the limit as $X\to\infty$ of $M_X(s;f)$, call it $M(s;f)$, exists and is holomorphic, by an application of Vitali's convergence theorem (see Section 5.21 of Titchmarsh \cite{Tit/1939} or Lemma C of \cite{Ka-Pe/resoI}) we obtain that the limit as $X\to\infty$ of $H_X(s;f)$ exists, is holomorphic and $\ll e^{\epsilon|t|}$ in any such substrip. As a consequence, $F(s;f)-M(s;f)$ has holomorphic continuation to $\si>-R$ and is $\ll e^{\epsilon|t|}$ there. Finally, Theorem 1 follows by further refining $M(s;f)$, similarly as in Section 2.5 on p.1421 of  \cite{Ka-Pe/2011a}.

\medskip
2.1.~{\bf Set up.} Let $F_X(s;f)$ be as above, $c>0$ be a constant (not necessarily the same at each occurrence and possibly depending on some parameters), $\epsilon>0$ be arbitrarily small and $\delta>0$ be sufficiently small. Moreover, for $N$ as above let $[N]=\{0,1,\dots,N\}$, $\emptyset \neq \A\subseteq [N]$, $|\A|$ be the cardinality of the set $\A$ and (with the obvious notation $\LL_{|\A}$)
\[
\bfw = \sum_{j=0}^N \kappa_jw_j, \quad \d \bfw = \prod_{j=0}^N \d w_j, \quad G(\bfw) = \prod_{j=0}^N \Gamma(w_j)z_j^{-w_j},
\]
\[
\bfw_{|\A} = \sum_{j\in\A} \kappa_jw_j, \quad \d \bfw_{|\A} = \prod_{j\in\A} \d w_j, \quad G(\bfw_{|\A}) = \prod_{j\in\A} \Gamma(w_j)z_j^{-w_j},
\]
\[
\frac{1}{2} <\eta<\frac{3}{4}, \quad \int_{\LL}\d\bfw=\int_{(-\eta)}\dots\int_{(-\eta)}\d\bfw \ \ \text{and analogously for} \ \int_{\LL_{|\A}}\d\bfw_{|\A},
\]
\[
I_X(s,\A)=\frac{1}{(2\pi i)^{|\A|}} \int_{\LL_{|\A}}F(s+\bfw_{|\A}) G(\bfw_{|\A}) \d \bfw_{|\A}.
\]
Let $K\geq 0$, $R\geq K+1$ be such that $\frac{d+1}{2}+dR\not\in\NN$, $\sigma>-R$ and $\rho>\frac{R+1}{\K}$, where $\K = \sum_{j=0}^N\kappa_j$. Then for $\Re w_j = \rho$ we have $\Re(s+\bfw) = \sigma+\K\rho>1$, hence by Mellin's transform we get
\begin{equation}
\label{2-1}
F_X(s;f) = \frac{1}{(2\pi i)^{N+1}} \int_{(\rho)} \dots \int_{(\rho)} F(s+\bfw)G(\bfw)\d \bfw.
\end{equation}
As in Lemma 2.1 of \cite{Ka-Pe/2011a} we want to shift the line of integration in \eqref{2-1} to $\LL$, but now we have to cross the possible pole of $F(s+\bfw)$ at $\bfw=1-s$, in addition to the poles of $G(\bfw)$ at $w_j=0$. This adds extra difficulties, hence we start with the following

\medskip
{\bf Lemma 2.1.} {\sl With the above notation, for $-R<\sigma<-R+\delta$ we have
\[
F_X(s;f) =  \sum_{\emptyset \neq \A \subseteq [N]} I_X(s,\A)  + H^{(1)}_X(s),
\]
where $H^{(1)}_X(s) \ll (1+|t|)^c$ as $|t|\to\infty$ with some $c>0$,  uniformly for $X\geq1$.}

\medskip
{\it Proof.} Suppose that $F(s)$ has the following Laurent expansion at $s=1$:
\[
F(s) = \sum_{k=1}^{m_F} \frac{\beta_k}{(s-1)^k} + F_0(s), \hskip1.5cm \text{$F_0(s)$ entire}.
\]
Let $-R<\si<-R+\delta$. We move the line of integration in the $w_0$-variable to the right (if necessary) to $\Re w_0 = \rho_0$ with $\rho_0>\frac{R+1}{\kappa_0}$, and all the other to the left to $\Re w_j=\epsilon$ with a small $\epsilon>0$. Since we don't cross any pole, choosing $\A=\{1,\dots,N\}$ we have
\begin{equation}
\label{2-2}
F_X(s;f) = \frac{1}{(2\pi i)^{N}} \int_{(\epsilon)} \dots \int_{(\epsilon)} \big(\frac{1}{2\pi i} \int_{(\rho_0)} F(s+\bfw)\Gamma(w_0) z_0^{-w_0} \d w_0 \big) G(\bfw_{|\A}) \d \bfw_{|\A}.
\end{equation}
Now we move the integration in the inner integral to $\Re w_0 = \epsilon$, hence such an integral equals
\begin{equation}
\label{2-3}
\begin{split}
R_X(s,w_1,\dots,w_N) &+ \frac{1}{2\pi i} \int_{(\epsilon)} F(s+\bfw)\Gamma(w_0) z_0^{-w_0} \d w_0, \\
R_X(s,w_1,\dots,w_N) &= \res_{w_0=\frac{1}{\kappa_0}(1-s-\bfw_{|\A})} F(s+\bfw) \Gamma(w_0) z_0^{-w_0}.
\end{split}
\end{equation}
The residual function $R_X(s,w_1,\dots,w_N)$ coincides with the residual function $R_N^{(1)}(s,\alpha)$ defined on p.333 of \cite{Ka-Pe/2005}, after the following changes in the latter function:
\begin{equation}
\label{2-4}
s\mapsto s+\bfw_{|\A}, \qquad d\mapsto 1/\kappa_0, \qquad N\mapsto X, \qquad \alpha\mapsto \alpha_0.
\end{equation}
Therefore, the computations on p.334 of \cite{Ka-Pe/2005} give (recall that the coefficients $\alpha_k$ of the Laurent expansion of $F(s)$ in \cite{Ka-Pe/2005} are now called $\beta_k$)
\[
R_X(s,w_1,\dots,w_N) = \sum_{k=1}^{m_F} \frac{\beta_k}{\kappa_0^k(k-1)!} \sum_{\nu=0}^{k-1} {k-1 \choose \nu} \Gamma^{(\nu)}\big(\frac{1-s-\bfw_{|\A}}{\kappa_0}\big) (-\log z_0)^{k-\nu-1} z_0^{\frac{s+\bfw_{|\A}-1}{\kappa_0}}.
\]
Hence, in view of \eqref{2-3}, the contribution of $R_X(s,w_1,\dots,w_N)$ to the integral in \eqref{2-2} equals
\begin{equation}
\label{2-5}
\begin{split}
\sum_{k=1}^{m_F} \frac{\beta_k}{\kappa_0^k(k-1)!} &\sum_{\nu=0}^{k-1} {k-1 \choose \nu} z_0^{\frac{s-1}{\kappa_0}} (-\log z_0)^{k-\nu-1}\times  \\
 &\times \frac{1}{(2\pi i)^{N}} \int_{(\epsilon)} \dots \int_{(\epsilon)} \prod_{j=1}^N\Gamma(w_j) (z_0^{-\frac{\kappa_j}{\kappa_0}}z_j)^{-w_j} \Gamma^{(\nu)}\big(\frac{1-s-\bfw_{|\A}}{\kappa_0}\big) \d \bfw_{|\A}.
 \end{split}
\end{equation}
Thanks to Stirling's formula and to the form of the $z_j$'s, recalling that $\kappa_j<\kappa_0$ we may shift the lines of integration in \eqref{2-5} to $-\infty$, thus getting that \eqref{2-5} equals
\begin{equation}
\label{2-6}
\sum_{k=1}^{m_F} \frac{\beta_k}{\kappa_0^k(k-1)!} \sum_{\nu=0}^{k-1} {k-1 \choose \nu} z_0^{\frac{s-1}{\kappa_0}} (-\log z_0)^{k-\nu-1} H_{\nu}(s,\bfz),
\end{equation}
where
\begin{equation}
\label{2-7}
\begin{split}
H_{\nu}(s,\bfz) &= \sum_{k_1=0}^\infty \cdots  \sum_{k_N=0}^\infty \frac{(-1)^{k_1+\cdots +k_N}}{k_1!\cdots k_N!} \Gamma^{(\nu)}\big(\frac{1-s+\sum_{j=1}^N\kappa_jk_j}{\kappa_0}\big) \prod_{j=1}^N (z_0^{-\frac{\kappa_j}{\kappa_0}}z_j)^{k_j} \\
&= \Gamma^{(\nu)}\big(\frac{1-s}{\kappa_0}\big) +  \widetilde{H}_{\nu}(s,\bfz),
\end{split}
\end{equation} 
say, where $\widetilde{H}_{\nu}(s,\bfz)$ denotes the above summation over $(k_1,\dots,k_N)\neq (0,\dots,0)$. Moreover, since $\kappa_j<\kappa_0$, from Stirling's formula we have that there exists a constant $\delta_1>0$ such that
\[
\Gamma^{(\nu)}\big(\frac{1-s+\sum_{j=1}^N\kappa_jk_j}{\kappa_0}\big) \ll \prod_{j=1}^N k_j^{(1-\delta_1)k_j}
\]
for $-R < \sigma < -R+\delta$ (the implicit constant may depend on $R$ and $\nu$). Hence from \eqref{2-7}
\[
\widetilde{H}_{\nu}(s,\bfz) \ll \prod_{j=1}^N\big(\sum_{k_j=0}^\infty \frac{c^{k_j}k_j^{(1-\delta_2)k_j}}{k_j!}\big) \ll 1
\]
for $-R < \sigma < -R+\delta$ and some $\delta_2>0$, uniformly for $X\geq 1$. As a consequence, the part of \eqref{2-6} coming from $ \widetilde{H}_{\nu}(s,\bfz)$, which we denote by $k_X(s)$, is uniformly bounded in $X$ for $-R<\sigma <-R+\delta$. Moreover,  the computations before (3.6) on p.334 of \cite{Ka-Pe/2005} and \eqref{2-4} show that the part of \eqref{2-6} coming from $ \Gamma^{(\nu)}\big(\frac{1-s}{\kappa_0}\big)$ equals
\begin{equation}
\label{2-8}
- \sum_{k=1}^{m_F}\frac{\beta_k}{(s-1)^k} + g_X(s),
\end{equation}
where $g_X(s)$ is also uniformly bounded in $X$ for $-R<\sigma <-R+\delta$. Therefore, gathering \eqref{2-2}, \eqref{2-3} and \eqref{2-5}-\eqref{2-8}, from the properties of $k_X(s)$ and $g_X(s)$ we get
\begin{equation}
\label{2-9}
F_X(s;f) =  - \sum_{k=1}^{m_F}\frac{\beta_k}{(s-1)^k}  + \frac{1}{(2\pi i)^{N+1}} \int_{(\epsilon)} \dots \int_{(\epsilon)} F(s+\bfw) G(\bfw) \d \bfw +h_X(s),
\end{equation}
where $h_X(s)$ is uniformly bounded for $X\geq1$ and $-R<\sigma <-R+\delta$.

\smallskip
Finally, shifting the lines to $-\eta$ we have to cross the poles of $G(\bfw)$ at $w_j=0$, hence we deal with the integral in \eqref{2-9} as in Lemma 2.1 of \cite{Ka-Pe/2011a}, thus getting
\[
\begin{split}
F_X(s;f) &=  F(s) - \sum_{k=1}^{m_F}\frac{\beta_k}{(s-1)^k}  + \sum_{\emptyset \neq \A \subseteq [N]}  \frac{1}{(2\pi i)^{|\A|}} \int_{\LL_{|\A}}F(s+\bfw_{|\A}) G(\bfw_{|\A}) \d \bfw_{|\A} +h_X(s) \\
&= \sum_{\emptyset \neq \A \subseteq [N]} I_X(s,\A) + H^{(1)}_X(s),
\end{split}
\]
say. Moreover, thanks to the properties of $h_X(s)$ and $F(s)$, $H^{(1)}_X(s)$ satisfies
\[
H^{(1)}_X(s) \ll (1+|t|)^c
\]
for some $c>0$, uniformly for $X\geq1$ and $-R<\sigma<-R+\delta$, and the lemma follows. \fine

\medskip
{\bf Remark 7.} Note that $H^{(1)}_X(s)$ is holomorphic for $-R<\sigma<-R+\delta$, but this information is not necessary to prove Theorem 1. Indeed, holomorphy of the relevant terms will follow in a simpler way at a later stage of the proof, as we pointed out in the sketch at the beginning of the section. The same remark applies to the other terms below, in the sense that they are holomorphic for $-R<\sigma<-R+\delta$, $\delta>0$ sufficiently small. Indeed, in general such terms are meromorphic with poles in a horizontal strip of finite height. We shall more or less implicitly use the latter property in the rest of the proof; we refer to Section 2 of \cite{Ka-Pe/2011a} for details. \fine

\medskip
Let $-R<\sigma<-R+\delta$. Writing $\bfx_{|\A} = s+\bfw_{|\A}$ and
\[
\widetilde{G}(\bfx_{|\A}) = \prod_{j=1}^r \Gamma(\lambda_j(1-\bfx_{|\A}) +\overline{\mu}_j) \Gamma(1-\lambda_j\bfx_{|\A} -\mu_j), \qquad S(\bfx_{|\A}) = \prod_{j=1}^r \sin\pi(\lambda_j\bfx_{|\A}+\mu_j),
\]
we apply, in the integrals $I_X(s,\A)$ , the functional equation to $F(s+\bfw_{|\A})$ and then the reflection formula to the resulting $\Gamma$-factors. Expanding the Dirichlet series of $\overline{F}(1-s-\bfw_{|\A})$ we obtain 
\begin{equation}
\label{2-10}
 I_X(s,\A) = \frac{\omega}{\pi^r}Q^{1-2s} \sum_{n=1}^\infty \frac{\overline{a(n)}}{n^{1-s}} \frac{1}{(2\pi i)^{|\A|}} \int_{\LL_{|\A}} \widetilde{G}(\bfx_{|\A}) S(\bfx_{|\A}) G(\bfw_{|\A}) \big(\frac{n}{Q^2}\big)^{\bfw_{|\A}} \d \bfw_{|\A};
\end{equation}
see (2.3) of \cite{Ka-Pe/2011a}. Since for $\bfw_{|\A}\in\LL_{|\A}$ we have $\Re(-\bfx_{|\A})>|\sigma|+\kappa_N\eta>0$, we may apply Stirling's formula to $\widetilde{G}(\bfx_{|\A})$. By a computation similar to the one leading to (2.4) of \cite{Ka-Pe/2011a} (recalling that here we do not assume that $\theta_F=0$ as in \cite{Ka-Pe/2011a}) for any integer $L>0$ we get
\begin{equation}
\label{2-11}
\widetilde{G}(\bfx_{|\A}) = B^{\bfy_{|\A}} \sum_{\ell=0}^L c_\ell \Gamma\big(\frac{d+1}{2} - d\bfy_{|\A}-\ell \big) + O\big(e^{-\frac{\pi}{2}d|\Im\bfx_{|\A}|} (1+|\Im\bfx_{|\A}|)^{-L+c}\big),
\end{equation}
where $\bfy_{|\A}=\bfx_{|\A}+i\theta_F$, $B = d^d/\beta$, $\beta=\prod_{j=1}^r\lambda_j^{2\lambda_j}$, $c_0\neq0$ and $c>0$. Turning to $S(\bfx_{|\A})$, by the same argument leading to (2.6) and (2.7) of \cite{Ka-Pe/2011a} we obtain
\begin{equation}
\label{2-12}
S(\bfx_{|\A}) = k_1e^{-i\frac{\pi}{2} d \bfx_{|\A}} + k_2e^{i\frac{\pi}{2} d \bfx_{|\A}} + O(e^{\frac{\pi}{2}(d-c)|\Im \bfx_{|\A}|})
\end{equation}
with constants $k_1,k_2\neq0$ and some $c>0$. Therefore, a computation based on \eqref{2-11} and \eqref{2-12} shows that
\begin{equation}
\label{2-13}
\widetilde{G}(\bfx_{|\A}) S(\bfx_{|\A}) = \big(k_1e^{-i\frac{\pi}{2} d \bfx_{|\A}} + k_2e^{i\frac{\pi}{2} d \bfx_{|\A}}\big) B^{\bfy_{|\A}} \sum_{\ell=0}^L c_\ell \Gamma\big(\frac{d+1}{2} - d\bfy_{|\A}-\ell \big) + R(\bfx_{|\A})
\end{equation}
with $R(\bfx_{|\A}) \ll 1$ provided $L>c$. Moreover, by Stirling's formula we have
\[
G(\bfx_{|\A}) \ll \prod_{j=0}^N (1+ |w_j|^{1/2+\eta})^{-1}
\]
uniformly for $X\geq1$ and hence, recalling that $\eta>1/2$, with the same uniformty we get
\begin{equation}
\label{2-14}
\frac{1}{(2\pi i)^{|\A|}} \int_{\LL_{|\A}} \big|R(\bfx_{|\A}) G(\bfx_{|\A}) \big(\frac{n}{Q^2}\big)^{\bfw_{|\A}} \d \bfw_{|\A} \big| \ll 1.
\end{equation}

\medskip
From \eqref{2-10}, \eqref{2-13} and \eqref{2-14} we finally obtain that for $-R<\sigma<-R+\delta$
\begin{equation}
\label{2-15}
I_X(s,\A) = e^{a_1s+b_1} \sum_{n=1}^\infty \frac{\overline{a(n)}}{n^{1-s}} \sum_{\ell=0}^L\big(e_\ell I_X(s,\A,n,\ell) + e'_\ell  J_X(s,\A,n,\ell)\big) + H_X^{(2)}(s,\A)
\end{equation}
where $a_1\in\RR$ and $b_1,e_\ell,e'_\ell$ are certain constants with $e_0e'_0\neq0$,
\begin{equation}
\label{2-16}
I_X(s,\A,n,\ell) = \frac{1}{(2\pi i)^{|\A|}} \int_{\LL_{|\A}} \Gamma\big(\frac{d+1}{2}-d\bfx_{|\A}-\ell -id\theta_F\big) e^{-i\frac{\pi}{2}d\bfx_{|\A}} G(\bfw_{|\A}) \big(\frac{n}{q}\big)^{\bfw_{|\A}} \d \bfw_{|\A},
\end{equation}
$q=Q^2/B$, $J_X(s,\A,n,\ell)$ is equal to $I_X(s,\A,n,\ell)$ with $e^{-i\frac{\pi}{2}d\bfx_{|\A}}$ replaced by $e^{i\frac{\pi}{2}d\bfx_{|\A}}$, and
\begin{equation}
\label{2-17}
\sum_{\emptyset\neq\A\subseteq [N]} H_X^{(2)}(s,\A)\ll 1
\end{equation}
as $|t|\to\infty$ uniformly for $X\geq1$.

\medskip
2.2.~{\bf Mellin transform.} Thanks to the factor $G(\bfw_{|\A})$, the integral \eqref{2-16} with $e^{- i\frac{\pi}{2}d\bfx_{|\A}}$ replaced by $e^{i\pi \Lambda\bfx_{|\A}}$, $|\Lambda|\leq d/2$, is meromorphic with poles in a horizontal strip of finite height. For $|\Lambda|<d/2$, the Mellin transform argument on p.1410 of \cite{Ka-Pe/2011a} shows that for $\frac{d+1}{2d}-\frac{\ell}{d} <\sigma< \frac{d+1}{2d}-\frac{\ell}{d} + \delta$
\begin{equation}
\label{2-18}
\begin{split}
&\frac{1}{(2\pi i)^{|\A|}} \int_{\LL_{|\A}} \Gamma\big(\frac{d+1}{2}-d\bfx_{|\A}-\ell - id\theta_F\big) e^{i\pi \Lambda\bfx_{|\A}} G(\bfw_{|\A}) \big(\frac{n}{q}\big)^{\bfw_{|\A}} \d \bfw_{|\A}  \\
&= f_\ell \int_0^\infty \exp(-e^{\pi i\Lambda/d}x^{1/d}) \prod_{j\in\A}\big(e^{-z_j(\frac{qx}{n})^{\kappa_j}}-1\big) x^{\frac{d+1}{2d} - \frac{\ell}{d} -s -1 -i\theta_F} \d x
\end{split}
\end{equation}
with a certain $f_\ell\neq0$. Indeed, under the above conditions all the integrals involved in the argument are absolutely convergent. Moreover, for $\Lambda=\mp d/2$ the integrals in \eqref{2-18} are also absolutely convergent for $s$ in the above range. Hence we let $\Lambda \to \mp d/2$ in \eqref{2-18}, thus getting similar expressions for $I_X(s,\A,n,\ell)$ and $J_X(s,\A,n,\ell)$. Next we sum such expressions over $\A$ and argue as on p.1411 of \cite{Ka-Pe/2011a} to get the analog of (2.17) and (2.19) of \cite{Ka-Pe/2011a}. We obtain
\begin{equation}
\label{2-19}
\sum_{\emptyset\neq\A\subseteq[N]} I_X(s,\A,n,\ell) =e^{a_2s} f_\ell \int_0^\infty  e^{ix^{1/d}} \big(e^{-\Psi_X(x,n)} e^{-2\pi if(\frac{qx}{n},\bfalpha)} -1\big) x^{\frac{d+1}{2d} - \frac{\ell}{d} -s -1 -i\theta_F} \d x
\end{equation}
for $\frac{d+1}{2d}-\frac{\ell}{d} <\sigma< \frac{d+1}{2d}-\frac{\ell}{d} + \delta$, where $a_2\in\RR$ and $\Psi_X(x,n) = \frac{1}{X} \sum_{j=0}^N \big(\frac{qx}{n}\big)^{\kappa_j}$, and similarly for $J_X(s,\A,n,\ell)$. Accordingly we write ($J_X(s,n,\ell)$ is as in \eqref{2-19} with $e^{-ix^{1/d}}$ in place of $e^{ix^{1/d}}$)
\begin{equation}
\label{2-20}
\sum_{\emptyset\neq\A\subseteq[N]} I_X(s,\A,n,\ell) = f_\ell I_X(s,n,\ell)  \ \ \text{and} \ \  \sum_{\emptyset\neq\A\subseteq[N]} J_X(s,\A,n,\ell) = f_\ell J_X(s,n,\ell).
\end{equation}

\medskip
The integral $J_X(s,n,\ell)$, without saddle point, is dealt with by the following

\medskip
{\bf Lemma 2.2.} {\sl $J_X(s,n,\ell)$ is meromorphic for $\sigma<\frac{d+1}{2d}+\delta$, and for every given $-R < \sigma< -R + \delta$ and $\epsilon>0$, $J_X(s,n,\ell) \ll e^{\epsilon|t|}$ as $|t|\to\infty$ uniformly for $n\in\NN$ and $X\geq 1$}.

\medskip
{\it Proof.} Let $\frac{d+1}{2d}-\frac{\ell}{d} <\sigma< \frac{d+1}{2d}-\frac{\ell}{d} + \delta$. We start as in the proof of Lemma 2.2 of \cite{Ka-Pe/2011a}, switching to the complex variable $z$ in place of the real variable $x$ and then replacing the path of integration by $z=\rho e^{-i\phi}$, $0<\rho<\infty$ and $\phi>0$ arbitrarily small. On the new path we have 
$|e^{iz^{1/d}}| = e^{-\rho^{1/d}\sin(\phi/d)}\ll e^{-c\phi\rho^{1/d}}$ and we split the integral as
\[
J_X(s,n,\ell) = \int_0^{e^{-i\phi}\infty} \hskip-.5cm \dots \d z = e^{-i\phi}\big(\int_0^1 \dots\d\rho + \int_1^{n^\epsilon} \hskip-.2cm\dots\d\rho +  \int_{n^\epsilon}^{\infty} \dots\d\rho\big) = I_1(s)+I_2(s) +I_3(s),
\]
say. Arguing similarly as in (2.21) of \cite{Ka-Pe/2011a}, using the above bound for $e^{iz^{1/d}}$ we see that $I_3(s)$ is entire and $\ll 1$ uniformly in $n$ and $X$ for $\sigma>-R$. Analogously, $I_2(s)$ is entire and satisfies
\[
I_2(s) \ll \int_1^{n^\epsilon} |e^{-\Psi_X(z,n)} e^{-2\pi if(\frac{qz}{n},\bfalpha)} -1| \rho^{\frac{d+1}{2d} +R -1} \d\rho \ll \frac{1}{n^{\kappa_N}} \int_1^{n^\epsilon} \rho^{\frac{d+1}{2d} +R -1+\kappa_N} \d\rho \ll 1
\]
uniformly in $n$ and $X$ for $\sigma>-R$, provided $\epsilon>0$ is sufficiently small.

\medskip
In order to deal with $I_1(s)$, for $0<\rho<1$ we expand the exponentials in such a way that given $W\geq L/d+1+\delta$ we find an integer $M\geq1$ with the property that
\[
e^{iz^{1/d}} \big(e^{-\Psi_X(z,n)} e^{-2\pi if(\frac{qz}{n},\bfalpha)} -1\big) = \sum_{m=1}^M \beta_m(n) \rho^{u_m} + E_M(z,n)
\]
with $|E_M(z,n)|\ll \rho^W$ uniformly in $n$ and $X$. Moreover, we also have that $\beta_m(n) \ll n^{-\kappa_N}$ and the exponents $u_m$ are of the form $\frac{k}{d} + \sum_{j=0}^N \ell_j\kappa_j$ with $k,\ell_j\geq0$, $ \sum_{j=0}^N \ell_j>0$ and $0<u_m<W$. Consequently we have
\[
I_1(s) =  \sum_{m=1}^M \frac{\beta_m(n) e^{-i\phi(\frac{d+1}{2d} - \frac{\ell}{d} -s -i\theta_F)}}{u_m + \frac{d+1}{2d} - \frac{\ell}{d} -s -i\theta_F} + h_{X,M}(s,n) = \Sigma_{X,M}(s,n) +  h_{X,M}(s,n),
\]
say, with
\[
h_{X,M}(s,n) = e^{-i\phi(\frac{d+1}{2d} - \frac{\ell}{d} -s -i\theta_F)} \int_0^1 E_M(z,n) \rho^{\frac{d+1}{2d} - \frac{\ell}{d} -s -1 -i\theta_F} \d\rho.
\]
Thanks to the choice of $W$, this integral is absolutely convergent for $\sigma< \frac{d+1}{2d} + \delta$. Hence in view of our choice of $\phi$, $h_{X,M}(s,n)$ is holomorphic and satisfies $h_{X,M}(s,n) \ll e^{\epsilon|t|}$ uniformly in $n$ and $X$ for $s$ in any finite vertical strip contained in  $\sigma< \frac{d+1}{2d} + \delta$. Clearly, $\Sigma_{X,M}(s,n)$ is meromorphic over $\CC$ with poles at $s= u_m + \frac{d+1}{2d} - \frac{\ell}{d} -s -i\theta_F$, $1\leq m\leq M$ and $0\leq \ell \leq L$, therefore we may choose $-R<\sigma<-R+\delta$ away from the poles. Moreover, thanks to the bound for $\beta_m(n)$, for such a $\sigma$ we have that $\Sigma_{X,M}(s,n) \ll e^{\epsilon|t|}$ uniformly in $n$ and $X$, and the lemma follows. \fine

\medskip
2.3.~{\bf Saddle point.} Here we follow closely the saddle point argument in Section 2.3 of \cite{Ka-Pe/2011a}, hence we only briefly outline the needed changes and refer to \cite{Ka-Pe/2011a} for details and notation (see also the Introduction for some notation). Let $\xi$ be sufficiently large, $x_0=x_0(\xi,\bfalpha)\in\RR$ be the critical point of $\Phi(z,\xi,\bfalpha)$ as in Lemma 2.3 of \cite{Ka-Pe/2011a} and
\begin{equation}
\label{2-21}
K_X(s,\xi) = \gamma x_0^{\frac{d+1}{2d}-s} \int_{-r}^r e^{-\Psi_X(z,\xi)+i\Phi(z,\xi,\bfalpha)}(1+\gamma\lambda)^{\frac{d+1}{2d}-s-1} \d\lambda,
\end{equation}
where $\gamma=1-i$, $z=x_0(1+\gamma\lambda)$ and $r\in(0,1)$ is given in (2.29) of \cite{Ka-Pe/2011a}. Clearly, $K_X(s,\xi)$ is an entire function since $\Re(1+\gamma\lambda)>0$. As in \cite{Ka-Pe/2011a} we show that, for $n$ sufficiently large, the main contribution to the meromorphic integral $I_X(s,n,\ell)$ comes from $K_X(s+\frac{\ell}{d}+i\theta_F,n)$.

\medskip
{\bf Lemma 2.3.} {\sl Let $n_0$ be sufficiently large. Then for $n\geq n_0$ we have
\[
I_X(s,n,\ell) = K_X(s+\frac{\ell}{d}+i\theta_F,n) + H_X^{(3)}(s,n,\ell),
\]
where $H_X^{(3)}(s,n,\ell)$ is meromorphic for $\sigma<\frac{d+1}{2d} + \delta$ and, for every given $-R < \sigma< -R + \delta$ and $\epsilon>0$, satisfies $H^{(3)}_X(\sigma+it,n,\ell) \ll e^{\epsilon|t|}$ as $|t|\to\infty$ uniformly for $n\geq n_0$ and $X\geq 1$. Moreover, $I_X(s,n,\ell)$ has the same properties of $H_X^{(3)}(s,n,\ell)$ for $n<n_0$.}

\medskip
{\it Proof.} This is the analog of Lemma 2.4 of \cite{Ka-Pe/2011a} and its proof is similar, provided we make the same variations we did in Lemma 2.2 above with respect to Lemma 2.2 in \cite{Ka-Pe/2011a}. In particular, we have to split the integral over the path $z=\rho e^{i\phi}$, $\phi$ arbitrarily small and $\rho>0$, into three parts, as for $J_X(s,n,\ell)$. We don't give details to keep the paper in a reasonable size. \fine

\medskip
From Lemma 2.1, \eqref{2-15}, \eqref{2-17}, \eqref{2-19}, \eqref{2-20} and Lemmas 2.2 and 2.3 we deduce that
\begin{equation}
\label{2-22}
F_X(s;f) = e^{as+b}  \sum_{\ell=0}^L g_\ell \sum_{n=n_0}^\infty \frac{\overline{a(n)}}{n^{1-s}} K_X(s+ \frac{\ell}{d} + i\theta_F,n) + H_X(s;f) = M_X(s;f) + H_X(s;f),
\end{equation}
say, where $n_0$ is sufficiently large, $a\in\RR$, $g_0\neq0$ and $b,g_1,\dots,g_L$ are constants. Moreover, for any given $-R<\sigma<-R+\delta$ and $\epsilon>0$, $H_X(\sigma+it;f) \ll e^{\epsilon|t|}$ uniformly for $X\geq1$.

\medskip
2.4.~{\bf Limit as $X\to\infty$.} We write EBV$(X)$ for ``entire and bounded on every vertical strip, depending on $X$''. Arguing exactly as in Lemma 2.5 of \cite{Ka-Pe/2011a} (the extra $i\theta_F$ we have here does not change the bounds in \cite{Ka-Pe/2011a}), recalling that $d\kappa_0>1$ and using \eqref{2-22} we have

\smallskip
$\bullet$ $H_X(s;f)$ is EBV$(X)$;

\smallskip
$\bullet$ for any given $\sigma>d\kappa_0$, $H_X(\sigma+it;f) \ll 1$ uniformly for $X\geq1$;

\smallskip
$\bullet$ for any given $-R<\sigma<-R+\delta$ and $\epsilon>0$, $H_X(\sigma+it;f) \ll e^{\epsilon|t|}$ uniformly for $X\geq1$.

\smallskip
\noindent
Hence, by an application of the Phragm\'en-Lindel\"of theorem to $\Gamma(2\epsilon s+c)H_X(s;f)$ ($c>0$ suitable) and the strip $\sigma_1\leq \sigma\leq \sigma_2$, where $-R<\sigma_1<-R+\delta$ and $\sigma_2>d\kappa_0$, we deduce that

\smallskip
(i) {\it for $\sigma_1,\sigma_2$ as above and every $\epsilon>0$, $H_X(s;f)$ is holomorphic in the strip $\sigma_1<\sigma<\sigma_2$ and satisfies $H_X(s;f)\ll e^{\epsilon|t|}$ as $|t|\to\infty$, uniformly for $X\geq1$.}

\smallskip
\noindent
Moreover, still arguing as in Lemma 2.5 of \cite{Ka-Pe/2011a}, we also have that

\smallskip
(ii) {\it the limit as $X\to\infty$ of $H_X(s;f)$ exists for every $s$ in the strip $d\kappa_0<\sigma<\sigma_2$.}

\medskip
Therefore, thanks to (i) and (ii), from Vitali's convergence theorem (see Section 5.21 of Titchmarsh \cite{Tit/1939} or Lemma C of \cite{Ka-Pe/resoI}) we deduce that
\[
H(s;f) = \lim_{X\to\infty} H_X(s;f)
\]
exists and is holomorphic in any substrip of $\sigma_1<\sigma<\sigma_2$, and satisfies $H(s;f)\ll e^{\epsilon|t|}$. Writing
\[
K(s,\xi) = \gamma x_0^{\frac{d+1}{2d}-s} \int_{-r}^r e^{i\Phi(z,\xi,\bfalpha)}(1+\gamma\lambda)^{\frac{d+1}{2d}-s-1} \d\lambda
\]
(notation is as in \eqref{2-21}), once again by the arguments in Section 2.5 of \cite{Ka-Pe/2011a}, applied to the term $M_X(s;f)$ in \eqref{2-22}, we have that 
\begin{equation}
\label{2-23}
M(s;f) = \lim_{X\to\infty} M_X(s;f) = e^{as+b}  \sum_{\ell=0}^L g_\ell \sum_{n=n_0}^\infty \frac{\overline{a(n)}}{n^{1-s}} K(s+ \frac{\ell}{d} + i\theta_F,n)
\end{equation}
exists and is holomorphic and bounded for $\sigma>d\kappa_0$. Thus, letting $X\to\infty$ in \eqref{2-22} we obtain
\[
F(s;f) = M(s;f) + H(s;f)
\]
where $H(s;f)$ is holomorphic for $\sigma>-K$ and satisfies $H(s;f)\ll e^{\epsilon|t|}$ as $|t|\to\infty$. 

\medskip
2.5.~{\bf Completion of the proof.} To complete the proof of Theorem 1 we follow the arguments in Section 2.5 of \cite{Ka-Pe/2011a}, but here we have to take into account many more terms in the expansions. Such terms will eventually contribute to the sum on the right hand side of \eqref{1-7}. Again, to keep the paper in a reasonable size, we outline the changes and refer to \cite{Ka-Pe/2011a} as much as possible. We recall that $x_0$ is as in \eqref{2-21}, $P=x_0^2\Phi''(x_0,\xi,\bfalpha)$ (see (2.29) of \cite{Ka-Pe/2011a}, where the same quantity is called $R$), $\gamma=1-i$ and $\kappa_0^*=\frac{\kappa_0}{d\kappa_0-1}$. Moreover, we denote by $Q(s,\xi)$ a finite sum of type
\begin{equation}
\label{2-24}
Q(s,\xi) = \sum_i P_{i}(s) r_{i}(\xi)
\end{equation}
with polynomials $P_{i}(s)$ and $r_{i}(\xi)=x_0^{a_{i}}/\xi^{b_{i}}$, $a_{i},b_{i}\geq0$. The analog of Lemma 2.7 of \cite{Ka-Pe/2011a} is

\medskip
{\bf Lemma 2.4.} {\sl Let $\xi_0$ be sufficiently large,  $\xi\geq \xi_0$, $0\leq \ell\leq L$, $M\geq 2$ be a given integer and $\epsilon>0$. Then there exist $h_j>0$ and $Q_{j,\ell}(s,\xi)$ as in \eqref{2-24}, $j=0,\dots,M(M+1)$, such that
\[
K(s+\frac{\ell}{d}+i\theta_F,\xi) =  \sum_{j=0}^{M(M+1)} \frac{\gamma h_j}{|P|^{(j+1)/2}} Q_{j,\ell}(s,\xi) x_0^{\frac{d+1}{2d} - \frac{\ell}{d} - s -i\theta_F} e\big(\frac{1}{2\pi}\Phi(x_0,\xi,\bfalpha)\big) + H_\ell^{(4)}(s,\xi)
\]
where $h_0=\sqrt{\pi}$, $Q_{0,\ell}(s,\xi)=1$ identically and $H_\ell^{(4)}(s,\xi)$ is entire. Moreover, for $s$ in any given vertical strip and $0\leq \ell\leq L$, $H_\ell^{(4)}(s,\xi)$ satisfies
\[
H_\ell^{(4)}(s,\xi) \ll e^{\epsilon|t|} \xi^{-d\kappa_0^*\sigma-(M-d)\kappa_0^*/2}\log^{3M+1}\xi,
\]
and the functions $r_{i}(\xi)$ inside $Q_{j,\ell}(s,\xi)$ satisfy $r_{i}(\xi)\ll \xi^{j\kappa_0^*/2}$.}

\medskip
{\it Proof.} Apart from switching from $R$ to $P$, we keep the notation of Lemma 2.7 of \cite{Ka-Pe/2011a} and write 
\begin{equation}
\label{2-25}
K(s+ \frac{\ell}{d} + i\theta_F,\xi) = \gamma x_0^{\frac{d+1}{2d} - \frac{\ell}{d} - s -i\theta_F} e^{i\Phi(x_0)} I(s), \quad  I(s) =  \int_{-r}^r e^{i (\Phi(z)-\Phi(x_0))} (1+\gamma\lambda)^{c(s,\ell)} \d\lambda
\end{equation}
where  $c(s,\ell) = \frac{d+1}{2d} - \frac{\ell}{d} - s-1 -i\theta_F$ and $r\ll \xi^{-\kappa_0^*/2}\log \xi$. We need the following expansions, based on (2.29), (2.30) and (2.41) of \cite{Ka-Pe/2011a} and valid for $-r\leq \lambda\leq r$. First we have
\[
\Phi(z)-\Phi(x_0) = \sum_{m=2}^M \frac{R_m}{m!} (\gamma\lambda)^m + O(\xi^{-(M-1)\kappa_0^*/2} \log^{M+1}\xi)
\]
(see p.1422 of \cite{Ka-Pe/2011a}). Since $R_2=P<0$ and  $| \sum_{m=3}^M (\gamma\lambda)^m\frac{R_m}{m!}|\ll \xi^{-\kappa_0^*/2} \log^3\xi$ we deduce that 
\[
\begin{split}
e&^{i (\Phi(z)-\Phi(x_0))} = e^{-|P|\lambda^2} \exp\big(i \sum_{m=3}^M \frac{R_m}{m!}(\gamma\lambda)^m\big)(1+O(\xi^{-(M-1)\kappa_0^*/2} \log^{M+1}\xi)) \\
&= \big(e^{-|P|\lambda^2} \sum_{k=0}^M \frac{1}{k!}\big(i \sum_{m=3}^M (\gamma\lambda)^m\frac{R_m}{m!}\big)^k +O(\xi^{-M\kappa_0^*/2} \log^{3M}\xi)\big) \big(1+O(\xi^{-(M-1)\kappa_0^*/2} \log^{M+1}\xi)\big) \\
&= \big(e^{-|P|\lambda^2} \big(1+ \sum_{h=3}^{M^2} R_h(\xi) \lambda^h\big) + O(\xi^{-M\kappa_0^*/2} \log^{3M}\xi) \big) \big(1+O(\xi^{-(M-1)\kappa_0^*/2} \log^{M+1}\xi)\big)
\end{split}
\]
where, in view of the definition of the $R_m$ on p.1421 of \cite{Ka-Pe/2011a}, the functions $R_h(\xi)$ are linear combinations of terms of type $x_0^a/\xi^b$. Moreover, for $s$ in any fixed vertical strip we have
\[
(1+\gamma\lambda)^{c(s,\ell)}= 1 + \sum_{k=1}^M P_k(s,\ell) \lambda^k +O(e^{\epsilon|t|} \xi^{-M\kappa_0^*/2} \log^{M}\xi)
\]
with certain polynomials $P_k(s,\ell)$. Hence, for $s$ in any vertical strip, the integrand in $I(s)$ equals
\[
e^{-|P|\lambda^2}\big(1 + \sum_{j=1}^{M(M+1)} Q_{j,\ell}(s,\xi)\lambda^j\big) + O(e^{\epsilon|t|} \xi^{-M\kappa_0^*/2}\log^{3M}\xi)
\]
where the $Q_{j,\ell}(s,\xi)$ are as in \eqref{2-24}. By (2.41) of \cite{Ka-Pe/2011a}, $|R_m|\ll \xi^{\kappa_0^*}$, hence $R_h(\xi)\ll \xi^{h\kappa_0^*/3}$ since the above summation over $m$ starts from $m=3$. Therefore, the functions $r_{i}(\xi)$ inside $Q_{j,\ell}(s,\xi)$ satisfy $r_{i}(\xi)\ll \xi^{j\kappa_0^*/2}$. Integrating over $[-r,r]$ and arguing similarly as on p.1422 of \cite{Ka-Pe/2011a} we get
\begin{equation}
\label{2-26}
I(s) = \sum_{j=0}^{M(M+1)} \frac{h_j}{|P|^{(j+1)/2}} Q_{j,\ell}(s,\xi) + O(e^{\epsilon|t|} \xi^{-(M+1)\kappa_0^*/2}\log^{3M+1}\xi)
\end{equation}
for $s$ in any vertical strip, where $h_0 = \sqrt{\pi}$, $h_j>0$ and $Q_{0,\ell}(s,\xi)=1$ identically. The lemma follows from \eqref{2-25} and \eqref{2-26}. \fine

\medskip
Now we are ready to conclude the proof. Recalling the notation before both \eqref{2-21} and Lemma 2.4, for $n\geq n_0$ we denote by $x_n$ the value of $x_0$ relative to $\xi=n$, i.e. $x_n=x_0(n,\bfalpha)$; analogously, $P_n$ is the value of $P$ relative to $\xi=n$. From \eqref{2-23} and Lemma 2.4 we obtain that
\[
M(s;f) = \gamma e^{as+b}  \sum_{\ell=0}^L g_\ell \sum_{j=0}^{M(M+1)} h_j \sum_{n=n_0}^\infty \frac{\overline{a(n)}}{n^{1-s}} \frac{Q_{j,\ell}(s,n)}{|P_n|^{(j+1)/2}}  x_n^{\frac{d+1}{2d} - \frac{\ell}{d} - s -i\theta_F} e\big(\frac{1}{2\pi}\Phi(x_n,n,\bfalpha)\big) + H^{(5)}(s;f)
\]
where, choosing $M=M(K)$ sufficiently large, $H^{(5)}(s;f)$ is holomorphic for $\sigma>-K$ and bounded by $e^{\epsilon|t|}$ for $s$ in any vertical strip inside $\sigma>-K$. Then we use Lemma 2.8 of \cite{Ka-Pe/2011a} and proceed analogously as on p.1424-1426. More precisely, we use the expansions of $1/\sqrt{|P_n|}$, $x_n^{\frac{d+1}{2d} - \frac{\ell}{d} - s -i\theta_F}$ and $e\big(\frac{1}{2\pi}\Phi(x_n,n,\bfalpha)\big)$ in (2.47), in the displayed equation after (2.47) and in (2.49) of \cite{Ka-Pe/2011a}, respectively. Thanks to the shape of the functions $r_i(\xi)$ in $Q_{j,\ell}(s,\xi)$ we also have
\[
Q_{j,\ell}(s,n) = n^{\theta_j}\sum_{\omega(j)} c_{\omega(j)}(s,\ell)n^{-\omega(j)}
\]
where, for each $j\geq1$, $0\leq\omega(j)\to\infty$ is a certain sequence, $0\leq \theta_j<j\kappa_0^*/2$ and $c_{\omega(j)}(s,\ell)$ is holomorphic and bounded by $e^{\epsilon|t|}$. Collecting sufficiently many terms in such expansions we get
\[
\frac{Q_{j,\ell}(s,n)}{|P_n|^{(j+1)/2}}  x_n^{\frac{d+1}{2d} - \frac{\ell}{d} - s -i\theta_F} e\big(\frac{1}{2\pi}\Phi(x_n,n,\bfalpha)\big) = \frac{e(f^*(n,\bfalpha))}{n^{d\kappa_0^*(s-\frac{d+1}{2d}+i\theta_F) +\frac{\kappa^*}{2}}} \sum_{\nu=0}^V \frac{c_{\nu,\ell,j}(s,\bfalpha)}{n^{\delta_{\nu,\ell,j}}} + H_{\ell,j}^{(6)}(s,n)
\]
where $c_{\nu,\ell,j}(s,\bfalpha)$ are entire with $c_{0,0,0}(s,\bfalpha)\neq0$, $\delta_{\nu,\ell,j}\geq0$ with $\delta_{0,0,0}=0$. Moreover, if $V=V(K)$ is sufficiently large, the sum over $n\geq n_0$ of the entire functions $\overline{a(n)}n^{s-1}H_{\ell,j}^{(6)}(s,n)$ is absolutely convergent for $\sigma>-K$ and bounded by $e^{\epsilon|t|}$ for $s$ in any vertical strip inside $\sigma>-K$. 

\medskip
Theorem 1 follows now summing the last equation over $n,j$ and $\ell$, since clearly
\[
s^* = 1 - s + d\kappa_0^*(s-\frac{d+1}{2d}+i\theta_F) +\frac{\kappa^*}{2}.
\]
The assertions in Remark 3 in the Introduction follow by an analysis of the above arguments.

\bigskip
\section{Proof of the other statements}

\smallskip
3.1.~{\bf Proof of Fact 1.} We proceed by induction on $j$, recalling that $\deg S\geq 1$; see after \eqref{1-8}. Clearly, for $j=1$ we have $\ell_1=\deg S_1$. Assuming that $\ell_{j-1}=\deg S_{j-1}$, thanks to \eqref{1-3} we have
\[
\text{lexp}((TS_{j-1}\cdots TS_1)^\flat) = \frac{\ell_{j-1}}{d\ell_{j-1}-1} <1,
\]
since $\deg S_{j-1}\geq 1$ and $d>2$. Hence applying $S_j$ we immediately have that $\ell_j=\deg S_j$. \fine

\bigskip
3.2.~{\bf Proof of Fact 2.} Suppose that $G(f_0)=H(f_1)$ and $G\neq H$, thus we may write
\[
\begin{split}
G &= TS_M\cdots TS_NTS_{N-1}\cdots TS_1 \\
H &= TS_M\cdots TS_NTS'_R\cdots TS'_1
\end{split}
\]
with, in particular, $S_{N-1}\neq S'_R$. Hence, since $\frak{G}$ is a group, we must have 
\[
S_{N-1}\cdots TS_1(f_0) = S'_R\cdots TS'_1(f_1).
\]
But $\deg(S_{N-1}^{-1}S'_R) \geq 1$, hence $f_0=S_1^{-1}T\cdots TS_{N-1}^{-1}S'_R\cdots TS'_1(f_1)$, all the operations being allowed by Fact 1. But this gives a contradiction, since $\text{lexp}(f_0)= \deg S_1^{-1}=\deg S_1\geq 1$ by Fact 1, and $\text{lexp}(f_0)\leq 1/d$ by definition of $\A(F)$. \fine

\bigskip
3.3.~{\bf Proof of Theorems 2, 3, 4 and 5.} We apply induction with respect to the weight $\omega(f)$. If $\omega(f)=0$ then $f(\xi) = f_0(\xi) + P(\xi)$ with $0\leq \ \text{lexp}(f_0) \leq 1/d$ and $P\in\ZZ[\xi]$. Hence $F(s;f)=F(s;f_0)$ and therefore the assertions of all theorems hold true in this case thanks to the results in \cite{Ka-Pe/resoI}, since $D(f)=1$ by definition in this case. In order to perform the inductive step, we suppose that $\omega(f)=M$ and assume the assertions of the theorems true for any $F\in\S^\sharp$ with degree $d\geq 1$ and any $f_1\in\A(F)$ (resp. $\A(F)\setminus\A_0(F)$, $\A_0(F)\setminus\A_{00}(F)$ and $\A_{00}(F)$) of weight $M-1$. Hence we have
\begin{equation}
\label{3-1}
f(\xi) = (f_1(\xi) + P(\xi))^*
\end{equation}
with some $P\in\ZZ[\xi]$ of degree $\geq 1$, lexp$(f_1+P)=\ell_M>1/d$ and $\omega(f_1)=M-1$. We shall use different arguments to show that $F(s;f)$ has the required properties according to the case at hand.

\medskip
To prove Theorem 2, suppose that $f,f_1\in\A(F)$ and satisfy \eqref{3-1}, and let $K\geq0$ be arbitrary. Then we apply Theorem 1, thus getting an expression of type \eqref{1-7} for $F(s;f)$. But in view of \eqref{3-1} and \eqref{1-6} we have
\begin{equation}
\label{3-2}
\overline{F}(s^*+\eta_j; f^*) = \overline{F}(s^*+\eta_j; f_1),
\end{equation}
hence the right hand side of \eqref{1-7} is meromorphic for $\sigma>-K$ by the inductive hypothesis. Moreover, still from Theorem 1 we have that
\[
W_j(s), G(s) \ll e^{\epsilon|t|}
\]
uniformly in every vertical strip contained in $\sigma>-K$. Theorem 2 follows from the inductive hypothesis, since $K$ is arbitrarily large. \fine

\medskip
To prove Theorem 3 we just observe that if $f\in\A(F)\setminus\A_0(F)$ then $f_1\in\A(F)\setminus\A_0(F)$ as well. Hence, since $F(s;f_1)$ is entire in this case, every $\overline{F}(s^*+\eta_j; f_1)$ is also entire, and arguing as before we obtain that $F(s;f)$ is entire as well, thus proving Theorem 3. \fine

\medskip
Suppose now that $f\in\A_0(F)\setminus\A_{00}(F)$. Then $f_1\in\A_0(F)\setminus\A_{00}(F)$ and by the inductive hypothesis, $F(s;f_1)$ has simple poles at $s_0 = \frac12 + \frac{1}{2dD(f_1)} -i\theta_F$ and on the half-line
\[
s=\sigma - i\theta_F \hskip1.5cm \sigma \leq \frac12 + \frac{1}{2dD(f_1)}.
\]
Writing $\kappa_0=$ lexp$(f)$, by \eqref{1-12} we have 
\begin{equation}
\label{3-3}
\kappa_0 = \frac{\ell_M}{(d\ell_M-1)}.
\end{equation}
Hence by \eqref{1-7}, \eqref{1-4} and \eqref{3-2}, recalling that $\theta_{\overline{F}} = -\theta_F$ we see that the poles of $F(s;f)$ are simple and lie on the half-line
\[
\frac{s+\frac{d\kappa_0}{2}-1+id\theta_F\kappa_0}{d\kappa_0-1} = \sigma +i\theta_F \hskip1.5cm \sigma \leq \frac12 + \frac{1}{2dD(f_1)}.
\]
Therefore, the polar half-line of $F(s;f)$ becomes
\[
s= (d\kappa_0-1)\sigma - \frac{d\kappa_0}{2} + 1 - i\theta_F = \sigma' - i\theta_F,
\]
say, where thanks to \eqref{3-3} we have
\[
\sigma' \leq (d\kappa_0-1)\big(\frac12 + \frac{1}{2dD(f_1)}\big) - \frac{d\kappa_0}{2} + 1 = \frac12 + \frac{d\kappa_0-1}{2dD(f_1)} = \frac12 + \frac{1}{2dD(f)}.
\]
Moreover, the initial point of such a half-line is a simple pole, and Theorem 4 follows. \fine

\medskip
Suppose finally that $f\in\A_{00}(F)$. Arguing as before we see that $F(s;f)$ has poles of order $\leq m_F$ on the half-line
\[
\frac{s+\frac{d\kappa_0}{2}-1+id\theta_F\kappa_0}{d\kappa_0-1} = \sigma -i\theta_F\big(\frac{(-1)^{M-1}}{D(f_1)}-1\big) \hskip1.5cm \sigma \leq \frac12 + \frac{1}{2D(f_1)},
\]
and a pole of order $m_F$ at its initial point. Hence
\[
s= (d\kappa_0-1)\sigma - \frac{d\kappa_0}{2} + 1 - i\theta_F\big(d\kappa_0 + \big(\frac{(-1)^{M-1}}{D(f_1)}-1\big)(d\kappa_0-1)\big) = \sigma' -i\theta_F D',
\]
say. But thanks to \eqref{3-3} we have
\[
D' = 1 + \frac{1}{d\ell_M-1} +  \big(\frac{(-1)^{M-1}}{D(f_1)}-1\big)  \frac{1}{d\ell_M-1}  = -  \big(\frac{(-1)^M}{D(f)}-1\big)
\]
and
\[
\sigma' \leq (d\kappa_0-1)\big(\frac12 + \frac{1}{2D(f_1)}\big) - \frac{d\kappa_0}{2} + 1 = \frac{d\kappa_0-1}{2D(f_1)} +\frac12 = \frac12 + \frac{1}{D(f)},
\]
thus proving Theorem 5. \fine

\bigskip
3.4.~{\bf Proof of Theorem 6.} Let $F\in\S_2^\sharp$. We start with the twist function $g(\xi) = k\xi^2 + \ell \xi + \beta \sqrt{\xi}$, where $k,\ell\in\ZZ$, $k>0$ and $\beta \in\RR$. Since clearly $g= S(f_0)$, where $f_0$ is the standard twist and $S$ is a shift of degree 2, we have that $g\in\A(F)$. Then we apply the operator $T$ and compute explicitly $T(g)^\flat=g^*\in\A(F)$; the function $f(\xi)$ in Theorem 6 will be closely related to $g^*(\xi)$.

\smallskip
To this end we recall the definition of $x_0$ and the displayed equation before (1.11) on p.1401 of \cite{Ka-Pe/2011a}, saying that 
\[
g^*(\xi) = \frac{1}{2\pi} \Phi(x_0,\xi,\bfalpha)^\flat,
\]
where $\Phi(z,\xi,\bfalpha)$ is as in \eqref{1-2} (with $g$ in place of $f$ in this case) and the real number $x_0=x_0(\xi)\geq 1$ is the unique solution (in a certain region) of the equation $\frac{\partial}{\partial z} \Phi(z,\xi,\bfalpha)=0$. Therefore, writing for simplicity $q$ for the conductor $q_F$ and $\Phi(z)$ for $ \Phi(z,\xi,\bfalpha)$, we first compute the critical point $x_0$ of the function
\begin{equation}
\label{3-4}
\Phi(z) = \big(1-2\pi \beta \sqrt{\frac{\lambda}{\xi}}\big) z^{1/2} - 2\pi\big(\frac{k\lambda^2}{\xi^2}z^2 + \frac{\ell\lambda}{\xi}z\big) \hskip1.5cm \lambda=\frac{q}{(4\pi)^2},
\end{equation}
which satisfies
\begin{equation}
\label{3-5}
\frac12\big(1-2\pi \beta\sqrt{\frac{\lambda}{\xi}}\big)x_0^{-1/2} = \frac{4\pi k \lambda^2}{\xi^2}x_0 + \frac{2\pi\ell\lambda}{\xi}.
\end{equation}
Putting $x_0^{1/2} = X_0$ we obtain the cubic equation
\[
X_0^3 + \frac{\ell \xi}{2k\lambda} X_0 - \frac{\xi^2\beta(\xi)}{8\pi k \lambda^2} = 0,
\]
where $\beta(\xi) = 1-2\pi \beta \sqrt{\lambda/\xi}$. Hence by Cardano's formulae we get
\begin{equation}
\label{3-6}
X_0 = a\xi^{2/3}\left\{\left(1+\sqrt{1+\frac{b}{\xi}}\right)^{1/3} + \left(1-\sqrt{1-\frac{b}{\xi}}\right)^{1/3} \right\}
\end{equation}
with
\begin{equation}
\label{3-7}
a=a(\xi)=\big(\frac{\beta(\xi)}{16\pi k \lambda^2}\big)^{1/3} \hskip2cm b =b(\xi)= \frac{32}{27} \frac{\pi^2\ell^3\lambda}{k\beta(\xi)^2}.
\end{equation}
In view of \eqref{3-4} and \eqref{3-5} we therefore have
\begin{equation}
\label{3-8}
\Phi(x_0) = \frac{3}{4}\beta(\xi)x_0^{1/2} - \frac{\pi \ell \lambda}{\xi} x_0 = \frac{3}{4}\beta(\xi)X_0 - \frac{\pi \ell \lambda}{\xi} X_0^2.
\end{equation}

\smallskip
In order to compute $g^*(\xi)$ we have to approximate the right hand side of \eqref{3-8} up to negative powers of $\xi$. To this end we first note that, since $b/\xi\to 0$ as $\xi\to\infty$, a simple computation gives
\[
\begin{split}
\big(1+\sqrt{1+\frac{b}{\xi}}\big)^{1/3} &= 2^{1/3} +O\big(\frac{1}{\xi}\big) \\
\big(1-\sqrt{1-\frac{b}{\xi}}\big)^{1/3} &= -\big(\frac{b}{2\xi}\big)^{1/3} \big(1+ O(\frac{1}{\xi})\big)
\end{split}
\]
as $\xi\to\infty$. Hence, observing that $a,b = O(1)$ as $\xi\to\infty$, from \eqref{3-6} we obtain
\[
\begin{split}
X_0 &=  a\xi^{2/3}\big(2^{1/3} - \big(\frac{b}{2\xi}\big)^{1/3}\big) + O(\xi^{-1/3}) \\
X_0^2 &= 2^{2/3} a^2 \xi^{4/3} - 2a^2 b^{1/3} \xi + O(\xi^{2/3}).
\end{split}
\]
Therefore from \eqref{3-8} we finally get
\begin{equation}
\label{3-9}
\begin{split}
\Phi(x_0) &= \frac{3}{4} \beta(\xi) a \xi^{2/3}\big(2^{1/3} - \big(\frac{b}{2\xi}\big)^{1/3}\big) - \pi\ell\lambda (2^{2/3} a^2 \xi^{1/3} - 2a^2 b^{1/3}) + O(\xi^{-1/3}) \\
&= 2\pi\big(A\xi^{2/3} + B\xi^{1/3} + C\xi^{1/6} + c_0 \big)+  O(\xi^{-1/6})
\end{split}
\end{equation}
as $\xi\to\infty$, where
\begin{equation}
\label{3-10}
A = \frac{3}{2(2kq^2)^{1/3}} \hskip1.5cm
B =  -\frac{3\ell}{4(2k^2q)^{1/3}} \hskip1.5cm
C = -\frac{\beta}{2(k^2q)^{1/6}}
\end{equation}
and $c_0$ is a certain constant. Indeed, writing
\[
\tilde{\beta} = 2\pi\beta\lambda^{1/2}, \quad \tilde{a} = (16\pi k \lambda^2)^{-1/3}, \quad \tilde{b} = \frac{32}{27}\frac{\pi^2\ell^3\lambda}{k}
\]
and recalling \eqref{3-7} and $\beta(\xi) = 1-\tilde{\beta}\xi^{-1/2}$, a computation shows that the relevant functions in the first line of \eqref{3-9} satisfy
\[
\begin{split}
\beta(\xi)a &=  \tilde{a} - \frac{4}{3} \tilde{\beta} \tilde{a}\xi^{-1/2} + O(\xi^{-1}), \\
\big(\frac{b}{\xi}\big)^{1/3} &=   \tilde{b}^{1/3}\xi^{-1/3} - \frac{2}{3} \tilde{\beta}\tilde{b}^{1/3} \xi^{-2/3} + O(\xi^{-4/3}), \\
a^2 &= \tilde{a}^2 + O(\xi^{-1/2}),  \\
a^2b^{1/3} &= c + O(\xi^{-1/2}),
\end{split}
\]
where $c$ is a certain constant. This proves the expansion of $\Phi(x_0)$ in the second line of \eqref{3-9}, apart from the explicit value of the constants $A,B$ and $C$ in \eqref{3-10}, since
\[
\Phi(x_0) = \big(\frac{3}{4} 2^{1/3} \tilde{a}\big) \xi^{2/3} - \big(\frac{3}{4} 2^{-1/3} \tilde{a}\tilde{b}^{1/3} + \pi\ell \lambda 2^{2/3} \tilde{a}^2\big) \xi^{1/3} -\big(2^{1/3}\tilde{\beta}\tilde{a}\big) \xi^{1/6} + c_0 + O(\xi^{-1/6}).
\]
Then a further computation (which we omit), involving also the value of $\lambda$ in \eqref{3-4}, proves \eqref{3-10}. Therefore
\[
g^*(\xi) = A\xi^{2/3} + B\xi^{1/3} + C\xi^{1/6} + c_0
\]
with $A,B,C$ as in \eqref{3-10} and a certain constant $c_0$, and $g^*\in\A_0(F) \Leftrightarrow |\beta|\in$ Spec$(F)$.

\smallskip
Note that since both $\ell$ and $\beta$ can be positive and negative, the two signs $-$ in \eqref{3-10} may be omitted.
Hence by the substitution $\alpha = \frac{\beta}{2(k^2q)^{1/6}}$ we see that Theorem 6 is proved in the case $k>0$, since the constant $c_0$ may clearly be omitted without affecting the statement. If $k<0$ we just recall \eqref{1-6}, and the result follows in its full generality. \fine

\newpage
\bigskip

\ifx\undefined\bysame{poly}.
\newcommand{\bysame}{\leavevmode\hbox to3em{\hrulefill}\ ,}
\fi

\vskip 1cm
\noindent
Jerzy Kaczorowski, Faculty of Mathematics and Computer Science, A.Mickiewicz University, 61-614 Pozna\'n, Poland and Institute of Mathematics of the Polish Academy of Sciences, 
00-956 Warsaw, Poland. e-mail: kjerzy@amu.edu.pl

\medskip
\noindent
Alberto Perelli, Dipartimento di Matematica, Universit\`a di Genova, via Dodecaneso 35, 16146 Genova, Italy. e-mail: perelli@dima.unige.it


\begin{thebibliography}{100} {\normalsize

\bibitem{Co-Gh/1993} J.B.Conrey, A.Ghosh - {\sl On the Selberg class of Dirichlet series: small degrees} - Duke Math. J. {\bf 72} (1993),  673--693.

\bibitem{Ka-Pe/1999a} J.Kaczorowski, A.Perelli - {\sl On the structure of the Selberg class, I: $0\leq d \leq 1$} - Acta Math. {\bf 182} (1999), 207--241.

\bibitem{Ka-Pe/2002} J.Kaczorowski, A.Perelli - {\sl On the structure of the Selberg class, V: $1<d<5/3$} - Invent. Math. {\bf 150} (2002), 485--516.

\bibitem{Ka-Pe/2005} J.Kaczorowski, A.Perelli - {\sl On the structure of the Selberg class, VI: non-linear twists} - Acta Arith. {\bf 116} (2005), 315--341.

\bibitem{Ka-Pe/2011a} J.Kaczorowski, A.Perelli - {\sl On the structure of the Selberg class, VII: $1<d<2$} - Annals of Math. {\bf 173} (2011), 1397--1441.

\bibitem{Ka-Pe/twist} J.Kaczorowski, A.Perelli - {\sl Twists, Euler products and a converse theorem for $L$-functions of degree 2} - To appear in Annali Scuola Normale Sup. Pisa;  arXiv:1207.2312.

\bibitem{Ka-Pe/resoI} J.Kaczorowski, A.Perelli - {\sl Twists and resonance of $L$-functions, I} - To appear in J. European Math. Soc.; arXiv:1304.4734

\bibitem{Tit/1939} E.C.Titchmarsh - {\sl The Theory of Functions} - second ed., Oxford U. P. 1939.

\bibitem{Zag/1985} D.Zagier - {\sl Modular points, modular curves, modular surfaces and modular forms} - Springer LN 1111 (1985), 225--248.

}
\end{thebibliography}
\end{document}